\documentclass{amsart}
\newtheorem{theorem}{Theorem}
\newtheorem{lemma}[theorem]{Lemma}
\newcommand{\cc}{\mathbb{C}}
\newcommand{\pp}{\mathbb{P}}

\title[Hartogs' theorem for projective spaces]
{Hartogs' theorem on separate holomorphicity for projective spaces}

\author{P. M. Gauthier}
\address{D\'epartement de Math\'ematiques et de Statistique,
Universit\'e de Montr\'eal, Case postale 6128, Succursale
Centre-ville, Montr\'eal, H3C~3J7, Canada.}
\email{gauthier@dms.umontreal.ca}

\author{E. S. Zeron}
\address{Depto. Matem\'aticas, CINVESTAV, Apartado
Postal 14-740, M\'exico D.F., 07000, M\'exico.}
\email{eszeron@math.cinvestav.mx}

\date{\today}
\thanks{Research supported by CRSNG(Canada),
Cinvestav(Mexico) and Conacyt(Mexico)}

\begin{document}
\begin{abstract}
If a mapping of several complex variables into projective space is
holomorphic in each pair of variables, then it is globally
holomorphic.
\end{abstract}

\maketitle
\section{Introduction}

Hartogs' Theorem \cite{Har} on separately holomorphic functions is
not valid for holomorphic mappings into a general complex manifold.
Some basic counterexamples involve the complex projective space
$\pp^1$. Let $g:\cc^2\to\pp^1$ be the mapping given by
$g(x,y)=[xy,x^2+y^2]$, with $g(0,0)=[0,1]$, in the homogeneous
coordinates of $\pp^1$. We have that $g$ is a mapping holomorphic in
each entry separately, but $g$ is not even continuous at the origin.
There are several works classifying those complex spaces for which
the Hartogs' Theorem on separately holomorphic mappings holds. A
short list includes \cite{LMH}, \cite{Sh1} and \cite{Sh2}.

No complex projective space $\pp^m$ satisfies the hypotheses presented
in \cite{LMH}, \cite{Sh1} or \cite{Sh2}. The main objective of this
work is to prove that a weak version of the Hartogs' Theorem holds for
mappings into complex projective spaces $\pp^m$.

By a coordinate $k$-plane $\Pi^k$, we understand any affine linear
subspace of $\cc^n$, with $n\geq{k}$, obtained by fixing $n-k$ of the
coordinates. The original Hartogs' theorem may be stated as follows: given
an open subset $\Omega\subset\cc^n$, and a function $f:\Omega\to\cc$ whose
restriction to each intersection $\Omega\cap\Pi^1$ is holomorphic, for
every coordinate $1$-plane $\Pi^1$, we have that $f$ is holomorphic on
$\Omega$. We may now present the main result of this paper.

\begin{theorem}[\textbf{Main}]\label{main}
Let $\Omega$ be an open domain in $\cc^n$, with $n\geq3$, and $\pp^m$
be a complex projective space. Given a mapping $f:\Omega\to\pp^m$ whose
restriction to each intersection $\Omega\cap\Pi^2$ is holomorphic, for
every coordinate $2$-plane $\Pi^2$; we have that $f$ is holomorphic on
$\Omega$.
\end{theorem}

This theorem will be proved in the next section. The Main Theorem does
not hold if we use continuity (or smoothness) instead of holomorphicity.
For example, consider the mapping $h:\cc^3\to\pp^1$ defined by
$h(x,y,z)=[xyz,|x|^3+|y|^3+|z|^3]$, with $h(0,0,0)=[0,1]$. The 
restriction of $h$ to every coordinate $2$-plane $\Pi^2$ is smooth, 
but $h$ is not even continuous at the origin.

\section{Proof of main theorem}

We begin by recalling the following theorem of Alexander, Taylor and
Ullman.

\begin{theorem}\cite[p.~340]{ATU}\label{atu}
Let $\Omega$ be an open subset of $\cc^{n+1}$. Let $r\geq1$ be a fixed
integer, and $A\subset\Omega$ be a subset such that the intersection
$A\cap\Pi^n$ with every coordinate $n$-plane $\Pi^n$ is either empty
or a subvariety of pure dimension $r$. Then, the set $A$ is closed
in $\Omega$.
\end{theorem}

\textbf{Remark.} The theorem in \cite[p.~340]{ATU} does not say
explicitly that the intersection $A\cap\Pi^n$ may be empty. However,
in their proof, they use the hypothesis that each set $W_j$ is a
$r$-dimensional subvariety of $\Omega_1$, recalling their own notation.
Looking carefully at their proof, we have that no set $W_j$ can be
empty, for the point $c_j$ is indeed contained in $W_j$. Hence, the
proof in \cite[p.~340]{ATU} works perfectly under the hypothesis that
each intersection $A\cap\Pi^n$ is either empty or a $r$-dimensional
subvariety.

\begin{lemma}\label{lem1}
Let $\Omega$ be an open domain in $\cc^{n+1}$, with $n\geq2$, and $\pp^m$
be a complex projective space. Given a function $f:\Omega\to\pp^m$ whose
restriction to each intersection $\Omega\cap\Pi^n$ is holomorphic, for
every coordinate $n$-plane $\Pi^n$, we have that $f$ is holomorphic on
$\Omega$.
\end{lemma}

\begin{proof} We shall prove that $f$ is holomorphic on a neighbourhood
of any fixed point $z_0\in\Omega$. Thus, we shall suppose from now on
that $\Omega$ is a polydisc, and that $f(z_0)=[1,0,0,\ldots]$ in the
homogeneous coordinates of $\pp^m$. Consider the open set
$U_0\subset\pp^m$ composed of the points $[1,\xi]$, for $\xi\in\cc^m$.
We only need to show that $f^{-1}(U_0)$ is an open neighbourhood of $z_0$
in $\Omega$. Then, a direct application of the original Hartogs' theorem
yields that $f$ is holomorphic on $f^{-1}(U_0)$, because $U_0$ is
obviously biholomorphic to $\cc^m$; and so we have that $f$ is
holomorphic on a neighbourhood of $z_0$.

Notice that the hyperplane at infinity $\pp^m\setminus{U_0}$ is the
set of all points $[0,\xi]$ with $\xi\neq0$. Define the set $E$ equal to
$f^{-1}(\pp^m\setminus{U_0})$. We assert that $E\cap\Pi^n$ is a subvariety
of $\Omega\cap\Pi^n$, for every coordinate $n$-plane $\Pi^n$. Let $z_1$ be
any given point in $E\cap\Pi^n$. We obviously have that $f(z_1)=[0,\xi]$
with $\xi\neq0$. We can suppose, without loss of generality, that its
second entry is different from zero. That is, $f(z_1)=[0,1,\eta]$ for
some $\eta\in\cc^{m-1}$. Consider the open set $U_1\subset\pp^m$ composed
of the points $[x,1,y]$, with $x\in\cc$ and $y\in\cc^{m-1}$. Notice that
$E\cap\Pi^n$ is a closed subset of $\Omega\cap\Pi^n$, and
$f^{-1}(U_1)\cap\Pi^n$ is an open neighbourhood of $z_1$ in
$\Omega\cap\Pi^n$ as well, because the restriction of $f$ to
$\Omega\cap\Pi^n$ is holomorphic. Besides, consider the holomorphic function
$\pi:U_1\to\cc$, defined by $\pi[x,1,y]=x$. We have that $\pi^{-1}(0)$ is
equal to $(\pp^m\setminus{U_0})\cap{U_1}$. Hence, the set $E\cap\Pi^n$ is
a subvariety of $\Omega\cap\Pi^n$ around $z_1$, because
$E\cap\Pi^n\cap{f}^{-1}(U_1)$ is equal to the inverse image
of zero under the holomorphic function $\pi\circ{f}|_{\Pi^n}$.

It is easy to deduce that $E\cap\Pi^n$ has only three possibilities:
it can either be empty, equal to $\Omega\cap\Pi^n$, or a subvariety of
$\Omega\cap\Pi^n$ with pure dimension $n-1$. Define $J$ to be the union
of all intersections $E\cap\Pi^n$ which are equal to $\Omega\cap\Pi^n$.
We assert that $J$ is a closed subset of $\Omega$. Firstly, we consider
only the coordinate $n$-planes whose first coordinate is constant, that
is, planes $\Pi^{n,0}$ of the form $\{x\}\times\cc^n$. Recall that
$\Omega$ is a polydisc, so it can be written as the product
$D_0\times\Delta_0$, with $D_0$ an open disc in $\cc$. Define
$\rho_0:\Omega\to{D_0}$ to be the projection on the first coordinate,
and $J_0$ to be the union of all sets $E\cap\Pi^{n,0}$ which are equal
to $\Omega\cap\Pi^{n,0}$. It is easy to deduce that
$$J_0\,=\,\bigcap_{y\in\Delta_0}\rho_0(E\cap(D_0\times\{y\}))\times\Delta_0.$$

Let $H^n$ be any coordinate $n$-plane which contains the line
$\cc\times\{y\}$. We know that $E\cap{H^n}$ is a closed subset of
$\Omega\cap{H^n}$. Whence, we also deduce that every $E\cap{H^n}$,
each $E\cap(D_0\times\{y\})$, and $J_0$ are all closed subsets of
$\Omega$. We may analyse, in the same way, the coordinate $n$-planes
$\Pi^{n,k}$ of the form $\cc^k\times\{x\}\times\cc^{n-k}$, for
$0\leq{k}\leq{n}$, and define $J_k$ to be union of all sets
$E\cap\Pi^{n,k}$ which are equal to $\Omega\cap\Pi^{n,k}$. At the end, we
deduce that each $J_k$ is closed in $\Omega$. Moreover, it is easy to
deduce that $J=\bigcup_kJ_k$ is also closed in $\Omega$.

Finally, consider $E^*:=E\setminus{J}$ and the open set
$\Omega^*:=\Omega\setminus{J}$ in $\cc^{n+1}$. Every intersection
$E^*\cap\Pi^n$ with a coordinate $n$-plane $\Pi^n$ is either empty or
an analytic set of pure dimension $n-1$. Then $E^*$ is a closed subset
of $\Omega^*$, after applying Theorem~\ref{atu}. The sets
$f^{-1}(U_0)$, $\Omega\setminus{E}$ and $\Omega^*\setminus{E^*}$ are
all equal, so $f$ is indeed holomorphic on the open neighbourhood
$f^{-1}(U_0)$ of $z_0$; as we wanted to prove.
\end{proof}

The proof of the Main Theorem is a direct application of Lemma~\ref{lem1}.

\begin{proof}\textbf{(Main Theorem).} Considering Lemma~\ref{lem1},
we may deduce that the restriction of $f$ to each intersection
$\Omega\cap\Pi^3$ is holomorphic, for every coordinate $3$-plane $\Pi^3$.
Proceeding by induction, we only need to apply Lemma~\ref{lem1} a finite
number of times in order to deduce that the restriction of $f$ to each
intersection $\Omega\cap\Pi^k$ is holomorphic, for every coordinate
$k$-plane $\Pi^k$, with $k\geq3$; and so $f$ if holomorphic on $\Omega$.
\end{proof}

\bibliographystyle{plain}

\begin{thebibliography}{10}

\bibitem{ATU} H.\ Alexander, B.\ A.\ Taylor, J.\ L.\ Ullman.
Areas of projections of analytic sets. \textit{Invent. Math.}
\textbf{16} (1972), pp.~335--341.

\bibitem{Har} F.\ Hartogs. Zur Theorie der analytischen Funktionen
mehrerer unabh\"angiger Ver\"anderlichen, insbesondere \"uber die
Darstellung derselben durch Reihen, welche nach Potenzen einer
Ver\"anderlichen fortschreiten. \textit{ Math. Ann.} \textbf{62}
(1906), no.~1, pp.~1--88.

\bibitem{LMH} Le Mau Hai, Nguyen Van Khue. Hartogs spaces, spaces
having the Forelli property and Hartogs holomorphic extension spaces.
\textit{Vietnam J. Math.} \textbf{33} (2005), no.~1, pp.~43--53.

\bibitem{Sh1} B.\ Shiffman. Hartogs theorems for separately
holomorphic mappings into complex spaces. \textit{C. R. Acad.
Sci. Paris S\'er. I Math.} \textbf{310} (1990), no.~3, pp.89--94.

\bibitem{Sh2} B.\ Shiffman. Separately meromorphic functions and
separately holomorphic mappings, pp.~191--198. \textit{Several
complex variables and complex geometry, Proc. Sympos. Pure Math.,
vol~52, Part 1.} Amer. Math. Soc., Providence, RI, 1991.

\end{thebibliography}

\end{document}